\input amstex
\input amsppt.sty
\magnification=\magstep1
\vsize=22.2truecm
\baselineskip=16truept
\NoBlackBoxes
\nologo
\pageno=1
\topmatter
\def\Z{\Bbb Z}
\def\N{\Bbb N}
\def\l{\left}
\def\r{\right}

\def\gs{\geqslant}
\def\bi{\binom}
\def\al{\alpha}
\def\l{\left}
\def\r{\right}
\def\bg{\bigg}
\def\({\bg(}
\def\[{\bg\lfloor}
\def\){\bg)}
\def\]{\bg\rfloor}

\def\f{\frac}
\def\Proof{\noindent{\it Proof}}

\def\Ack{\medskip\noindent {\bf Acknowledgments}}

\hbox{J. Number Theory 130(2010), no.\,12, 2701--2706.}
\medskip
\title  On 2-adic orders of some binomial sums \endtitle
\author
Hao Pan and Zhi-Wei Sun
\endauthor
\address
Department of Mathematics, Nanjing University,
Nanjing 210093, People's Republic of China
\endaddress
\email{haopan79\@yahoo.com.cn,\ zwsun\@nju.edu.cn}\endemail
\abstract We prove that for any nonnegative integers $n$ and $r$ the
binomial sum
$$
\sum_{k=-n}^n\binom{2n}{n-k}k^{2r}
$$
is divisible by $2^{2n-\min\{\alpha(n),\alpha(r)\}}$, where
$\alpha(n)$ denotes the number of 1s in the binary expansion of $n$.
This confirms a recent conjecture of Guo and Zeng [J. Number Theory,
{\bf 130}(2010), 172--186].
\endabstract
\thanks 2010 {\it Mathematics Subject Classification}.\,Primary 11B65;
Secondary 05A10, 05A19, 11A07, 11S99.
\newline\indent The first author was supported by the National Natural Science
Foundation for Youths in China (grant 10901078).
\newline\indent The second author was the corresponding author; he was
supported by the National Natural Science
Foundation (grant 10871087) and the Overseas Cooperation Fund (grant
10928101) of China.
\endthanks
\endtopmatter
\document
\TagsOnRight

In 1976 Shapiro [3] introduced the Catalan triangle $(\f kn\bi{2n}{n-k})_{n\gs k\gs1}$
and determined the sum of entries in the $n$th row; namely, he showed that
$$\sum_{k=1}^n k\bi{2n}{n-k}=\f n2\bi{2n}n.$$

Let $n,r\in\N=\{0,1,2,\ldots\}$.  Recently, Guo and Zeng [1]
proved that
$$\f2{n^2\bi{2n}n}\sum_{k=1}^n\binom{2n}{n-k}k^{2r+1}$$
is an odd integer if $n,r\in\Z^+=\{1,2,3,\ldots\}$. They also conjectured that the binomial sum
$$
F(n,r)=\sum_{k=-n}^n\binom{2n}{n-k}k^{2r}\tag1.1
$$
is divisible by $2^{2n-\min\{\alpha(n),\alpha(r)\}}$, where
$\alpha(n)$ denotes the number of 1s in the binary expansion of $n$.
Note that if $n,r\in\Z^+$ then
$F(n,r)=2\sum_{k=1}^n\bi{2n}{n-k}k^{2r}$. Actually  the conjecture
was motivated by Guo and Zeng's following observations:
$$\align \sum_{k=1}^n\bi{2n}{n-k}k^2=&2^{2n-2}n,
\\\sum_{k=1}^n\bi{2n}{n-k}k^4=&2^{2n-3}n(3n-1),
\\\sum_{k=1}^n\bi{2n}{n-k}k^6=&2^{2n-4}n(15n^2-15n+4),
\\\sum_{k=1}^n\bi{2n}{n-k}k^8=&2^{2n-5}n(105n^3-210n^2+147n-34).
\endalign$$

In this paper we shall confirm the sophisticated conjecture of Guo and
Zeng.

For an integer $n$ and a prime $p$, the $p$-adic order of $n$
at $p$ is given by
$$
\nu_p(n)=\sup\{v\in\N:\,p^v\mid n\}.
$$

Now we state our main result.

 \proclaim{Theorem 1.1} For any $n,r\in\N$ we have
$$
\nu_2(F(n,r))\geqslant2n-\min\{\alpha(n),\alpha(r)\},\tag 1.2
$$
where $F(n,r)$ is given by $(1.1)$.
\endproclaim

Note that (1.2) can be split into two inequalities:
$$\nu_2(F(n,r))\geqslant2n-\alpha(n)\tag1.3$$
and $$ \nu_2(F(n,r))\geqslant2n-\alpha(r).\tag1.4$$

In Sections 2 and 3 we will show (1.3) and (1.4) respectively.

\heading{2. Proof of (1.3)}\endheading

Let $p$ be any prime. A useful theorem of Legendre (see, e.g., [2,
pp. 22--24]) asserts that for any $n\in\N$ we have
$$\nu_p(n!)=\sum_{i=1}^\infty\l\lfloor\frac n{p^i}\r\rfloor=\f{n-\al_p(n)}{p-1},$$
 where $\al_p(n)$ is the sum of the digits of $n$ in the expansion of $n$ in base $p$.
 In particular, $\nu_2(n!)=n-\al(n)$ for all $n=0,1,2,\ldots$.

\proclaim{Lemma 2.1} {\rm (i)} For any $n\in\Z^+$ we have
$$\nu_2(n)-1=\alpha(n-1)-\alpha(n).\tag2.1$$

{\rm (ii)} Let $s>t\gs0$ be integers. Then
$$\nu_2\bigg(\binom{s}{t}\bigg)\gs\alpha(t)-\alpha(s)+1.\tag2.2$$
\endproclaim
\Proof. (i) In view of Legendre's theorem, for any positive integer
$n$ we have
$$\nu_2(n)=\nu_2(n!)-\nu_2((n-1)!)=n-\al(n)-(n-1-\al(n-1))=\al(n-1)-\al(n)+1.$$
This proves (2.1).

(ii) With the help of Legendre's theorem,
$$\align
\nu_2\bigg(\binom{s}{t}\bigg)=&\nu_2(s!)-\nu_2(t!)-\nu_2((s-t)!)
\\=&s-\al(s)-(t-\al(t))-(s-t-\al(s-t))
\\=&\al(t)-\al(s)+\al(s-t)
\\\gs&\al(t)-\al(s)+1 \ (\text{since}\ s-t\gs1).
\endalign$$
So (2.2) holds. \qed

\proclaim{Lemma 2.2} For $n,r\in\Z^+$ we have
$$
F(n,r)=n^2F(n,r-1)-2n(2n-1)F(n-1,r-1).\tag 2.3
$$
\endproclaim
\Proof. Since
$$
(n^2-k^2)\binom{2n}{n-k}=2n(2n-1)\binom{2n-2}{n-1-k},
$$
we have
$$
\sum_{k=-n}^n\binom{2n}{n-k}k^{2r}=n^2\sum_{k=-n}^n\binom{2n}{n-k}k^{2r-2}
-2n(2n-1)\sum_{k=-n+1}^{n-1}\binom{2n-2}{n-1-k}k^{2r-2},
$$
which gives (2.3). \qed

\medskip
\noindent{\it Proof of (1.3)}.  We use induction on $n+r$. Clearly
(1.3) holds trivially when $n=0$ or $r=0$.

Now let $n,r\in\Z^+$ and assume (1.3) for any smaller value of
$n+r$. By (2.1), (2.3) and the induction hypothesis, we have
$$
\align
\nu_2(F(n,r))\geqslant&\min\{\nu_2(n^2F(n,r-1)),\nu_2(2n(2n-1)F(n-1,r-1))\}\\
=&\min\{2\nu_2(n)+\nu_2(F(n,r-1)),1+\nu_2(n)+\nu_2(F(n-1,r-1))\}\\
\geqslant&\min\{2\nu_2(n)+2n-\alpha(n),1+\nu_2(n)+2(n-1)-\alpha(n-1)\}\\
=&2n-\alpha(n).
\endalign
$$
This concludes the induction step. \qed

\heading{3. Proof of (1.4)}\endheading

\proclaim{Lemma 3.1} For  $n,r\in\Z^+$ we have
$$
\align F(n,r)=&4F(n-1,r)-\sum_{i=0}^{r-1}\binom{2r}{2i}F(n,i)-
2(2n-1)\sum_{i=0}^{r-1}\binom{2r}{2i+1}F(n-1,i)\\
&+n\sum_{i=0}^{r-1}\binom{2r}{2i+1}F(n,i)+2\sum_{i=0}^{r-1}\binom{2r}{2i}F(n-1,i).\tag
3.1\endalign
$$
\endproclaim
\Proof. Let $n\in\N$ and $r\in\Z^+$. We want to prove (3.1) with $n$
in it replaced by $n+1$.

Clearly
$$
\align
F(n,r)=&\sum_{k=-n-1}^{n-1}\binom{2n}{n-1-k}(k+1)^{2r}\\
=&\sum_{k=-n-1}^{n}\bigg(\binom{2n+1}{n-k}-\binom{2n}{n-k}\bigg)\sum_{j=0}^{2r}\binom{2r}{j}k^j\\
=&\sum_{j=0}^{2r}\binom{2r}{j}\sum_{k=-n-1}^{n}\binom{2n+1}{n-k}k^{j}
-\sum_{i=0}^{r}\binom{2r}{2i}\sum_{k=-n}^{n}\binom{2n}{n-k}k^{2i};
\tag 3.2
\endalign
$$
in the last step we use the fact that if $j$ is odd then
$$
\align
\sum_{k=-n}^n\binom{2n}{n-k}k^{j}=&\frac12\bigg(\sum_{k=-n}^n\binom{2n}{n-k}k^{j}
+\sum_{k=-n}^n\binom{2n}{n+k}(-k)^{j}\bigg)\\
=&\sum_{k=-n}^n\binom{2n}{n-k}(k^{j}+(-k)^{j})=0.
\endalign$$
When $j$ is even, we have
$$
\align
2\sum_{k=-n-1}^n\binom{2n+1}{n-k}k^{j}=&\sum_{k=-n-1}^{n+1}\bigg(\binom{2n+1}{n-k}+\binom{2n+1}{n+k}\bigg)k^{j}\\
=&\sum_{k=-n-1}^{n+1}\binom{2n+2}{n+1-k}k^{j}=F\left(n+1,\frac
j2\right).\tag 3.3
\endalign
$$
If $j$ is odd, then
$$
\align
&(n+1)\sum_{k=-n-1}^n\binom{2n+1}{n-k}k^{j-1}+\sum_{k=-n-1}^n\binom{2n+1}{n-k}k^{j}\\
=&\sum_{k=-n-1}^n(n+1+k)\binom{2n+1}{n-k}k^{j-1}=(2n+1)\sum_{k=-n}^n\binom{2n}{n-k}k^{j-1},
\endalign
$$
i.e.,
$$
\align
\sum_{k=-n-1}^n\binom{2n+1}{n-k}k^{j}=&(2n+1)\sum_{k=-n}^n\binom{2n}{n-k}k^{j-1}
-(n+1)\sum_{k=-n-1}^n\binom{2n+1}{n-k}k^{j-1}\\
=&(2n+1)F\left(n,\frac{j-1}2\right)-\frac{n+1}{2}F\left(n+1,\frac{j-1}2\right),\tag
3.4
\endalign
$$
where we use (3.3) in the last step. Combining (3.2)-(3.4), we get
$$
\align F(n,r)=&\frac{1}{2}\sum_{i=0}^{r}\binom{2r}{2i}F(n+1,i)+
(2n+1)\sum_{i=0}^{r-1}\binom{2r}{2i+1}F(n,i)\\
&-\frac{n+1}{2}\sum_{i=0}^{r-1}\binom{2r}{2i+1}F(n+1,i)-\sum_{i=0}^{r}\binom{2r}{2i}F(n,i),
\endalign
$$
which yields the desired result. \qed

\medskip
\noindent{\it Proof of (1.4)}. We still use induction on $n+r$.
There is nothing to do if $n=0$ or $r=0$. Assume that $n,r\geqslant
1$ and (1.4) holds for any smaller value of $n+r$. In view of Lemma
3.1,  $\nu_2(F(n,r))$ is not smaller than the minimum of the
following numbers:
$$
\align&2+\nu_2(F(n-1,r)),\
\min_{0\leqslant i<r}\nu_2\bigg(\binom{2r}{2i}F(n,i)\bigg),\
\min_{0\leqslant i<r}\nu_2\bigg(n\binom{2r}{2i+1}F(n,i)\bigg)\\
&1+\min_{0\leqslant i<r}\nu_2\bigg(\binom{2r}{2i+1}F(n-1,i)\bigg),\
1+\min_{0\leqslant i<r}\nu_2\bigg(\binom{2r}{2i}F(n-1,i)\bigg).
\endalign
$$
By the induction hypothesis and Lemma 2.1(ii), we have
$\nu_2(F(n-1,r))\geqslant 2n-2-\alpha(r)$, and also
$$
\nu_2\bigg(\binom{2r}{2i}F(n,i)\bigg)\geqslant
2n-\alpha(i)+\alpha(2i)-\alpha(2r)+1=2n-\alpha(r)+1,
$$
$$
\nu_2\bigg(\binom{2r}{2i+1}F(n-1,i)\bigg)\geqslant
2n-2-\alpha(i)+\alpha(2i+1)-\alpha(2r)+1=2n-\alpha(r),
$$
$$
\nu_2\bigg(n\binom{2r}{2i+1}F(n,i)\bigg)\geqslant2n-\alpha(i)+\alpha(2i+1)-\alpha(2r)+1=2n-\alpha(r)+2,
$$
and
$$\nu_2\bigg(\binom{2r}{2i}F(n-1,i)\bigg)\geqslant2n-2-\alpha(i)+\alpha(2i)-\alpha(2r)+1=2n-\alpha(r)-1.
$$
Thus (1.4) follows.\qed

\medskip
\Ack. The authors wish to thank the referees for their helpful comments.

\enddocument